\documentclass[a4paper]{article}
\usepackage{typearea}
\areaset[1cm]{158mm}{234mm}
\usepackage{xcolor}
\usepackage{graphicx}
\usepackage{psfrag}
\usepackage{caption}
\usepackage{subcaption}

\usepackage{fancybox}
\usepackage{fancyvrb}

\usepackage{amsmath}
\usepackage{amsthm}
\usepackage{amssymb}

\usepackage{latexsym,amssymb}
\usepackage{url}
\usepackage{amsfonts}
\usepackage{hhtensor}

\usepackage{indentfirst}

\def\RR{\mathbb{R}}

\def\HH{\mathbb{H}}
\def\LL{\mathbb{L}}
\def\MM{\mathbb{M}}

\def\TT{\mathbb{T}}

\def\G{{\cal G}}
\def\D{{\cal D}}
\def\J{{\cal J}}
\def\I{{\cal I}}

\def\M{{\cal M}}

\def\L{{\cal L}}
\def\F{{\cal F}}

\def\dfrac{\displaystyle\frac}
\def\tens{\mathbf}
\def\tens{\boldsymbol}
\def\media_#1{\mathchoice{\mathop{\hskip2pt{-}\hskip-10.3pt\int\hskip-4pt}%
    \nolimits_{#1}}%
    {\mathop{-\hskip -8.8pt\int}\nolimits_{#1}}{}{}}

\let\ffi\varphi

\def\div{\mathrm{div}}

\def\Im{\mathrm{Im}}
\def\dint{\displaystyle\int}

\theoremstyle{plain}
\newtheorem{theorem}{Theorem}

\newtheorem{lemma}[theorem]{Lemma}

\newtheorem{proposition}[theorem]{Proposition}

\theoremstyle{definition}
\newtheorem{definition}[theorem]{Definition}

\newtheorem{remark}[theorem]{Remark}

\title{Stress formulation and duality approach in periodic homogenization}
\author{Cristian Barbarosie$^1$\thanks{\texttt{cabarbarosie@fc.ul.pt}} \and Anca-Maria Toader$^2$\thanks{\texttt{atoader@fc.ul.pt}, 
both authors are from DM-FCUL and CMAFcIO Centro de Matem\'atica, Aplica\c c\~oes Fundamentais e Investiga\c c\~ao Operacional, 
Departamento de Matem\'atica da Faculdade de Ci\^encias, 
Universidade de Lisboa, 1749-016 Lisboa, Portugal.
}}
\date{\today}

\begin{document}
	\maketitle
	
	\begin{abstract}

This paper describes several different variational formulations of the so-called ``cellular problem''
which is a system of partial differential equations arising in the theory of homogenization,
subject to periodicity boundary conditions.
These variational formulations of the cellular problems have as main unknown the
displacement, the stress or the strain, respectively.
For each of these three cases, an equivalent minimization problem is introduced.
The variational formulation in stress proves to have a distinguished role and it gives rise
to two dual formulations, one in displacement-stress and another one in strain-stress.
The corresponding Lagrangians may be used in numerical optimization for developing algorithms
based on alternated directions, of Uzawa type.
		



  

\end{abstract}


\section*{Introduction}
\label{sec:intro}

In engineering, it is common practice to approach linearized elasticity
problems through formulations in displacement and in stress.
Using the stress as unknown is convenient because, above a critical stress value,
the material's behaviour gets out from the linearized elasticity framework.
For an optimal design {\sl the principal stresses form naturally a flow chart between the load points},
similar to the velocity of a fluid.
This observation follows from analysing the famous Michell trusses \cite{Michell1904}.
The concept of fully stressed design \cite[Section 9.1]{HG1992} is also aligned with this observation,
although it is still an empirical principle.

This is why, for optimality purposes, we need to understand the stress formulation
in relation with the formulation in displacement and
in strain, as dual formulations to the primal formulation in stress,
according to the Legendre-Fenchel duality theory.
The two lagrangians thus obtained have a saddle-point;
this unifying theory was developed by Ciarlet in \cite{Ciarlet-G-K2011} where a new duality approach to
elasticity is stated and proven.
In the present work we extend this duality approach to periodic problems
which is the natural framework of the
cellular problems that define the homogenized elastic tensor according to the theory of homogenization.

We believe that the obtained Lagrangians may play a role in numerical methods for optimizing
the properties of the homogenized elastic tensor, namely in an alternated directions algorithm
of Uzawa type, as in \cite[subsection 8.1]{BenGoLie2005}.

\section{Periodic functions}
\label{sec:period}

We begin by describing the general periodicity notion that we have been using since \cite{BT2010}.

Consider a lattice of vectors in $ \RR^3 $
(an additive group $ \G $ generated by three independent vectors $ \vec g_1, \vec g_2, \vec g_3 $).
Define the parallelipiped  $ Y \subset \mathbb{R}^3 $\,:
$$
Y = \{ s_1 \vec g_1 + s_1 \vec g_2 + s_1 \vec g_3 : s_1, s_2, s_3 \in [0, 1]\}.
$$
$Y$ is called ``periodicity cell''.

A function $\ffi $ is said to be $\G$-periodic (or $Y$-periodic) if it is invariant to translations
with vectors in $ \G $.

\begin{remark} \label{rem:int-tau-Y}
  Let $ \psi $ be a $ \G $-periodic function; let $ \tau $ be a translation of $ \mathbb R^3 $
  (not necessarily belonging to $ \G $). Then
  $$ \int_{\tau(Y)} \psi = \int_Y \psi \,. $$  
\end{remark}
  
In a periodic microstructure, the rigidity is a fourth order tensor field
$ \tens C : \mathbb{R}^3 \mapsto \mathbb{R}^{81} $ which is $\G$-periodic, 
that is, it varies according to a periodic pattern.
Denote by $ \tens D $ the compliance tensor, that is, the inverse tensor of $ \tens C $.

According to the homogenization theory, the macroscopic behaviour of a body with periodic structure
is described by a constant homogenized elastic tensor $ \tens C^H $.
It is possible to define $ \tens C^H $ in terms of the solutions of the cellular problems
which are PDEs subject to periodicity conditions:


\begin{equation}  \label{eq:cell-pb}
\left\{\begin{matrix}
\div \bigl(\tens C\, \tens e(u_A)\bigr)  =  0 \hbox{ in } \RR^3 \\
u_A(y)  =  A y + \ffi_{ A}(y), 
\quad \hbox{ with }  \ffi_A \quad \G-\hbox{periodic} ,
\end{matrix}\right.
\end{equation}
where $A$ is a given macroscopic strain (a $ 3 \times 3 $ symmetric matrix) and $\tens e$
represents the corresponding microscopic strain, that is, the symmetric part of the gradient
of $ u_A $.
In the sequel, for historical reasons, we shall use also the notation $ \nabla\!{}_s $ for the
symmetric part of the gradient, when the intention is to focus on the displacement $ u_A $.

The solution $ u_A $ depends linearly on the matrix $A$ and verifies 
\begin{equation} \label{eq:A}
A = \media_Y  \tens e  (u_A) = \frac{1}{|Y|} \int_Y  \tens e (u_A) \,.
\end{equation}

\noindent The macroscopic stress associated to $ u_A $ is defined by
\begin{equation}  \label{eq:S}
S = \media_Y \tens C\, \tens e  (u_A) = \frac{1}{|Y|} \int_Y \tens C\, \tens e (u_A) 
\end{equation}
and consequently depends linearly on $A$.
The symbol $ \displaystyle\media_Y $ denotes the average value of the integrand
on the periodicity cell $Y$.
	
The homogenized elastic tensor is then defined, for all macroscopic strains $A$,
through the equality $\tens C^H A = S $, that is,
\begin{equation} \label{CHA}
\tens C^H A = \media_Y \tens C\, \tens e  (u_A) .
\end{equation}

\noindent $\tens C^H$ can be also defined, equivalently, in terms of energy products\,:
\begin{equation} \label{CHAB}
\langle \tens C^H A, B \rangle = \media_Y \langle \tens C\, \tens e (u_A), \tens e (u_B)\rangle,
\end{equation}
where $A$ and $B$ are any two symmetric matrices (strains).
See formula (\ref{eq:CH A,B}) below for a better comprehension.

Denote by $ \tens D^H $ the homogenized compliance tensor, that is, $ \tens D^H = (\tens C^H)^{-1} $
is the inverse of the homogenized elastic tensor.

\section{Fundamental ingredients}
\label{sec:Korn-Green}

We begin this section by presenting some notation used thoughout the paper.

Denote by $ L^2_\# (\RR^3,\RR^3) $ the space of vector fields in $ L^2_{\mbox{\small loc}} (\RR^3,\RR^3) $
which are $ \G $-periodic, endowed with the $ L^2 (Y,\RR^3) $ norm.
Denote by $ H^1_\# (\RR^3,\RR^3) $ the space of vector fields in $ H^1_{\mbox{\small loc}} (\RR^3,\RR^3) $
which are $ \G $-periodic, endowed with the $ H^1 (Y,\RR^3) $ norm.
Equivalently, $ H^1_\# (\RR^3,\RR^3) $ can be viewed as the completion
in the norm of $ H^1(Y,\RR^3) $ of the space of functions in $ C^\infty(\RR^3, \RR^3) $ which are 
$\G$-periodic.
We also define the space $ H^1_{\#0} (\RR^3,\RR^3) = \Bigl\{\ffi\in H^1_\# (\RR^3,\RR^3) \mid
\displaystyle \media_Y \ffi =0 \Bigr\} $.
Denote by $ H^{-1}_\# (\RR^3, \RR^3) $ the dual space of $ H^1_\# (\RR^3, \RR^3) $.

For the sake of abreviation we shall use the notations $ L^2_\# $ for $ L^2_\# (\RR^3,\RR^3) $,
$H^1_\# $ for $H^1_\# (\RR^3,\RR^3)$ and $ H^{-1}_\# $ for $ H^{-1}_\# (\RR^3, \RR^3) $.

Denote by $ \LL^2_\# $ the space of matrix fields in $ L^2_{\mbox{\small loc}} (\RR^3,\RR^9) $
which are $ \G $-periodic, endowed with the $ L^2 (Y,\RR^9) $ norm;
we add the subscript $s$ for symmetric matrices\,: $ \LL^2_{\# s} $.
Denote by $ \HH^1_\# $ the space of matrix fields in $ H^1_{\mbox{\small loc}} (\RR^3,\RR^9) $
which are $ \G $-periodic, endowed with the $ H^1 (Y,\RR^9) $ norm;
we add the subscript $s$ for symmetric matrices\,: $ \HH^1_{\# s} $.

Denote by $\HH_{\# s} (\div )$ the space

$$
\HH_{\# s} (\div ) = \bigl\{ \tens\mu \in \LL^2_{\# s} \mid \div \tens\mu \in  L^2_\# \bigr\}\,,
$$

\noindent where, with Einstein's repeated index notation, the operator
$ \div : L^2_{\mbox{\small loc}} (\RR^3,\RR^9) \to H^{-1}(\RR^3, \RR^3) $ is defined by
$ \div \tens\mu = (\mu_{ij,j})_{1\le i\le 3} $.

Consider $ LP $ the space of linear plus periodic displacements defined in $\RR^3$\,:
$$
\begin{matrix}
LP = \bigl\{ u : \RR^3 \to \RR^3 \mid u(y)  = Ay + \ffi (y), \\
 A \in \M_3^s(\RR), \ffi \in H^1_\# \bigr\},
\end{matrix} 
$$
where $ \M_3^s(\RR) $ denotes the space of symmetric $ 3\times 3 $ matrices with real elements.
$ LP $ is a Hilbert
space since it is a direct sum between a finite dimensional space and $ H^1_\# $.

For an arbitrarily fixed strain matrix $ A\in  \M_3^s(\RR) $, denote by $ LP(A) $ the set of linear plus periodic 
displacements having the linear part $ Ay $ \,:
$$
\begin{matrix}
  LP(A)= \bigl\{ u : \RR^3 \to \RR^3 \mid u(y) = Ay + \ffi (y), \ffi \in H^1_\# \bigr\}
  = \{ Ay \} + H^1_\# .
\end{matrix} $$
Thus the last equation in (\ref{eq:cell-pb}) is equivalent to $ u_A\in LP(A) $.
Note that $ LP(0) = H^1_\#  $ is a closed subspace of $ LP $ and for
a given strain matrix $ A\in \M_3^s(\RR) $ the set $ LP(A) $ is a translation of $ LP(0) $. 
Moreover,
$$ LP = \displaystyle\bigcup_{A\in \M_3^s(\RR)} LP(A). $$ 

Define the space of linear plus periodic functions having zero average by\,: 
\begin{equation}
\label{eq:LP_0}
LP_0 :=  \{ v \in LP \mid \media_Y v dy =0 \}.
\end{equation}

For a given $ 3\times 3 $ symmetric matrix $ S \in \M_3^s(\RR) $, define the space of stress fields 
\begin{equation} \label{eq:def-TT-S}
\TT(S) = \{ \tens\mu \in \HH_{\# s}(\div) \mid \div \tens\mu = 0, \quad \media_Y\tens\mu\, dy = S  \}.
\end{equation}

Consider the space $ \TT $ defined as
\begin{equation} \label{eq:def-TT}
\TT : = \bigcup_{S\in \M_3^s(\RR)} \TT(S) =
\{ \tens\mu \in \HH_{\# s} (\div ) \mid \div \tens\mu = 0 \hbox{ in } L^2_\# \} . 
\end{equation}

\begin{remark}\label{rem:LPA-avrg}
  For a given symmetric matrix $A$, a function $ v\in LP $ belongs to $ LP(A) $ if
  and only if $ \displaystyle \media_Y \tens e(v) = A $.
  See \cite[Lemma 2]{BT2010}.
  Similarly, by definition (\ref{eq:def-TT-S}), a matrix field $ \tens\mu\in \TT $ belongs to
  $ \TT_S $ if and only if $ \displaystyle \media_Y \tens\mu = S $.
\end{remark}

By $ \TT(0) $ we denote the space of stresses that have zero mean.
Define the orthogonal complement of $ \TT(0)\, $:
\begin{equation} \label{eq:def-TT-0-ort}
\TT(0)^\perp  : = \bigl\{ \tens e \in \LL^2_{\# s}  \mid \dint_Y \langle \tens e, \tens\mu \rangle =0,
\forall \tens\mu \in  \TT(0) \bigr\} . 
\end{equation}

In the remaining of this section we present some ingredients needed to prove the results in
sections \ref{sec:var-form}, \ref{sec:cell-pb-min} and \ref{sec:dual prob}.
Among the stated results are the Korn inequality, the Green's formula,
extensions of the Donati representation Theorem, all adapted to the periodic framework.
The proofs of these results may be found in \cite{BT2022-A} jointly with other fundamental results
adapted to periodic context in relation to the homogenization theory.

We begin by stating the Korn inequality in the periodic context, which stays behind 
the isomorphisms in Theorems \ref{thrm:isomorphism} and \ref{thrm:LP-isomorphic-TT_0-perp}.

\begin{theorem}
\label{thrm:Korn}
The Korn inequality below holds for a positive constant $C$ and for all $v$ in $H^1_{\# }$
\begin{equation}
\label{eq:Korn_ineq}
\|v\|_{H^1_\#} \le C \bigl( \| v \|_{\LL^2_\# } + \| \nabla\!{}_s v \|_{\LL^2_\# }\bigr)
\end{equation}
\end{theorem}

\begin{theorem} \label{thrm:Green}
The Green's formula 
\begin{equation}  \label{eq:Green}
\dint_Y \langle \tens\mu , \nabla\!{}_s v \rangle\, dy + \dint_Y (\div \tens\mu) \cdot v\, dy = 0
\end{equation}
holds for all $\tens\mu \in \HH_{\# s} (\div) $ and all $v \in H^1_\# (\RR^3,\RR^3)$.		
\end{theorem}

The proof of Green's formula above is based on the notion of trace and its properties in the context
of periodic functions, see \cite[Section 4]{BT2022-A}.

Theorem \ref{thrm:prod-avrg} below is a direct consequence of Green's formula.
But it also has far-reaching connections with the compensated compactness and
the div-curl Lemma introduced by L. Tartar in \cite{Tartar1979}, 
as pointed out in \cite[Section 5]{BT2022-A}.

\begin{theorem}   \label{thrm:prod-avrg}
Given an arbitrary $ v \in LP $ and an arbitrary $ \tens\sigma \in \TT $, the following equality holds 
\begin{equation}  \label{eq:prod-avrg-sym}
\media_Y \langle \nabla\!{}_s v, \tens\sigma \rangle\, dy =
\bigl\langle \media_Y \nabla\!{}_s v\, dy, \media_Y \tens\sigma\, dy \bigr\rangle .
\end{equation}
\end{theorem}

The above Theorem \ref{thrm:prod-avrg} may be stated in the following form that is due to P.\ Suquet, see \cite{Suquet1987}\,:

 {\it If $ v\in LP(A) $ and $ \tens\sigma \in \TT(S) $, 
\begin{equation} \label{eq:AS}
  \media_Y \langle \nabla\!{}_s v, \tens\sigma \rangle\, dy = \langle A, S \rangle .
\end{equation}
}

Indeed, for a given macroscopic stress $S$, the space $ \TT(S) $ is defined as the set of
  stress matrix fields in $ \TT $ whose mean is $S$.
  Similarly, for given macroscopic strain $A$, the space $ LP(A) $ is equal to the set of
  strain matrix fields in $ LP $ whose mean is $A$.
  The notation $ LP $ haz been used by the authors in \cite[Lemma 2]{BT2010}, the only difference is that in the present paper
the symmetry of the matrix $A$ is imposed in the very definition of the space $ LP $

For two arbitrary symmetric matrices (strains) $ A$ and $ B$, the equality (\ref{eq:AS}) implies that the homogenized tensor $ \tens C^H $ 
can be defined, according to (\ref{CHAB}), by
\begin{equation} \label{eq:CH A,B}
\langle \tens C^H  A ,  B \rangle =
\bigl\langle \media_Y \tens C\, \tens e  (w_{ A}),  B \bigr\rangle =
\media_Y \langle \tens C\, \tens e  (w_{ A}),  B \rangle =
\media_Y \langle \tens C\, \tens e  (w_{ A}), \tens e  (w_{ B})\rangle .
\end{equation}
The following three results, Theorems \ref{thrm:Donati.1}, \ref{thrm:Donati.2} and \ref{thrm:Donati.3}, 
are extensions of Donati's representation Theorem and generalise to the periodic 
framework Theorems 4.2 and 4.3 in \cite{Amrouche2006}. Their proofs and some specific remarks may be found in \cite{BT2022-A} .

\begin{theorem} \label{thrm:Donati.1}
Consider $\tens e \in \LL^2_{\# s}$. 
Then there exists $ v \in H^1_\# $ such that $ \tens e = \nabla\!{}_s v $ in $ \LL^2_{\# s} $
if and only if 
\begin{equation}
\dint_Y \langle \tens e , \tens s\rangle\, dy =0, \hbox{ for all } \tens s \in \LL^2_{\# s} \hbox{ such that } \div\, \tens s =0 \hbox{ in }
H^{-1}_\# .
\label{eq:Donati.1}
\end{equation}
In this case $v$ is unique up to an additive constant vector.
If we add the zero-average hypothesis, $ v \in H^1_{\#0} $, then $v$ is unique.
\end{theorem}

\begin{theorem} \label{thrm:Donati.2}
Let $ \tens e \in \LL^2_{\# s} $. Then 
\begin{equation}
\dint_Y \langle \tens e, \tens\mu \rangle =0 \hbox{ for all } \tens\mu \in \TT(0)
\label{eq:Donati.2}
\end{equation}
if and only if there exists $ v \in LP $ such that $ \tens e= \nabla\!{}_s v $.
In this case $v$ is unique up to an additive constant vector.
\end{theorem}

While in Theorem \ref{thrm:Donati.1} the mean of $ \tens e $ is zero, 
in Theorem \ref{thrm:Donati.2} the mean of $ \tens\mu $ is zero.
We can let both $ \tens e $ and $ \tens\mu $ to have non-zero average and state the following
Theorem \ref{thrm:Donati.3} below that
join together Theorems \ref{thrm:Donati.1} and \ref{thrm:Donati.2}.
Moreover the following result containts both Theorem \ref{thrm:prod-avrg} and its reciprocal.

\begin{theorem} \label{thrm:Donati.3}
Consider $\tens e \in \LL^2_{\# s}$. 
Then there exists $ v \in LP $ such that $ \tens e = \nabla\!{}_s v $ in $ \LL^2_{\# s} $
if and only if 
\begin{equation}
  \media_Y \langle \tens e , \tens s\rangle\, dy =
  \langle \media_Y \nabla\!{}_s v\, dy, \media_Y \tens s\, dy \rangle
  \hbox{ for all } \tens s \in \LL^2_{\# s} \hbox{ such that } \div\, \tens s = 0 .
\label{eq:Donati.3}
\end{equation}
\end{theorem}


The following two results allow to 
conclude the equivalence between the formulation in displacement and the formulation in strain in the next Section \ref{sec:cell-pb-min} .

\begin{theorem}[$H^1_{\# 0}$ is isomorphic to $\TT^\perp$]
\label{thrm:isomorphism}
The operator $\nabla\!{}_s : H^1_{\# 0} \mapsto \LL^2_{\# s}$ is an isomorphism from the space 
$$
H^1_{\# 0} : = \{ v \in  H^1_\# (\RR^3,\RR^3) : \media_Y v\, dy =0  \}
$$
onto $\Im \nabla\!{}_s $. Consequently $\Im \nabla\!{}_s $ is closed in $\LL^2_{\# s}$, more precisely
$\Im \nabla\!{}_s = \TT^\perp $.

Moreover, the dual operator of $\nabla\!{}_s: H^1_{\# 0} \mapsto \LL^2_{\# s}$ is 
$-\div : \LL^2_{\# s} \mapsto H^{-1}_\#$.
\end{theorem}

The above result is a consequence of Theorems 22, 26 and 27 in \cite{BT2022-A}
and is an ingredient of the proof of Theorem \ref{thrm:Donati.1}.
Jointly with the next theorem, it gives a more profound knowledge on the space of strain fields
$ \TT^\perp $.

The following result states an isomorphism between the space $LP_0$ and the orthogonal complement
of $\TT(0)$ by mean of the inverse operator of $\nabla\!{}_s$.
It is an important ingredient of the formulation in strain obtained in Theorem \ref{thrm:form-strain}
from the next Section \ref{sec:cell-pb-min} .

\begin{theorem}[$LP_0$ is isomorphic to $\TT(0)^\perp$]
\label{thrm:LP-isomorphic-TT_0-perp}
For each $\tens e \in \TT(0)^\perp $, denote by $\F(\tens e)$ the unique element in $LP_0$ that satisfies
$ \nabla\!{}_s  \F(\tens e) = \tens e$ (according to Theorem \ref{thrm:Donati.2}).
Then the mapping $\F : \TT(0)^\perp  \mapsto LP_0 $ defines an isomorphism between the Hilbert spaces
$\TT(0)^\perp$ and $LP_0$.
\end{theorem}

The above two isomorphisms allow us to introduce two variational formulations in strain,
see Propositions \ref{thrm:var-form-strain-1} and \ref{thrm:var-form-strain-2} in the next section.

\section{Variational formulations of the cellular problem}
\label{sec:var-form}

In this section we state and prove several equivalent variational formulations
where the main unknown is the displacement, the strain and the stress, respectively.
This approach gives a novel perspective in the field of periodic homogenization showing that
each of these quantities -- displacement, strain or stress -- completely characterize the
properties of the homogenized material.
In these formulations, either a macroscopic strain $A$ or a macroscopic stress $S$ is given as data;
the equivalence between these formulations is guaranteed by the subjacent relation
$ S = \tens C^H A $, equivalent to $ A = \tens D^H S $.

\begin{proposition} \label{prop:var-form-1}
The variational formulation of the cellular problem (\ref{eq:cell-pb}), when the macroscopic strain
$ A \in \M^s_3(\RR) $ is given, has the following form\,:

\begin{equation} \label{eq:var-form-displ-macrostrain}
\left\{
\begin{matrix} \mbox{find } u_A \in LP(A) \mbox{ such that } \\
  \displaystyle\int_{Y} \bigl\langle \tens C\, \tens e(u_A) ,\tens e(v)\bigr\rangle\, dy =0
  \ \ \forall v \in H^1_\# \,. 
\end{matrix} \right. 
\end{equation}
\end{proposition}

\begin{proof}
  The second line in (\ref{eq:cell-pb}) says that $ u_A \in LP(A) $.
  To prove that the first line in (\ref{eq:cell-pb}) implies the second line in
  (\ref{eq:var-form-displ-macrostrain}) it suffices to apply Green's formula (Theorem \ref{thrm:Green}).
  The difficult part is to prove that the second line in (\ref{eq:var-form-displ-macrostrain})
  implies the first line in (\ref{eq:cell-pb}).

  Suppose $ u_A $ satisfies (\ref{eq:var-form-displ-macrostrain}).
  By choosing a test function $ v\in \D(Y) $, the space of $ C^\infty $ functions with compact support,
  we conclude that $ \div \bigl( \tens C\, \tens e(u_A) \bigr) = 0 $ almost everywhere
  for the Lebesgue measure in $Y$.
  However, this does {\em not\,} imply $ \div \bigl( \tens C\, \tens e(u_A) \bigr) = 0 $
  almost everywhere in $ \mathbb R^3 $; $ \div \bigl( \tens C\, \tens e(u_A) \bigr) $
  could be a measure (or a distribution) concentrated on the interface between different translations
  of $Y$ (with vectors in the periodicity group $ \G $); see Remark \ref{rem:u''+2}.

  However, the periodic character of $ \tens C $, $ \tens e(u_A) $ and $ \tens e(v) $
  ensures that $ u_A $ satisfies not only (\ref{eq:var-form-displ-macrostrain}) but also
  $$ \int_{\tau(Y)} \bigl\langle \tens C\, \tens e(u_A) ,\tens e(v)\bigr\rangle\, dy =0
  \ \ \forall v \in H^1_\# \,, $$
  for any translation $ \tau $ of $ \mathbb R^3 $ (see Remark \ref{rem:int-tau-Y}).
  We stress that $ \tau $ does not have to belong to the periodicity group.
  Thus, $ \div \bigl( \tens C\, \tens e(u_A) \bigr) = 0 $ almost everywhere
  for the Lebesgue measure in $ \tau(Y) $.
  Since the translation $ \tau $ is arbitray, we conclude that 
  $ \div \bigl( \tens C\, \tens e(u_A) \bigr) = 0 $ almost everywhere in $ \RR^3 $.

\end{proof}

\begin{figure}[ht] \centering
  \psfrag{x0}{${}_0$}
  \psfrag{x1}{${}_1$}
  \psfrag{x2}{${}_2$}
  \psfrag{x3}{${}_3$}
  \psfrag{x4}{${}_4$}
  \includegraphics[width=100mm]{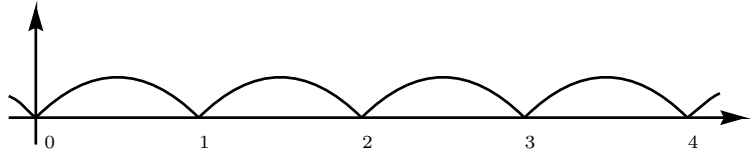}
  \caption{A periodic function with $ u''=-2 $ in $ {]}0,1{[} $}
  \label{fig:bumps}
\end{figure}

\begin{remark} \label{rem:u''+2}
  A function $ u\in H^1_\# $ may satisfy the equality $ \div \bigl(\tens C\, \tens e(u)\bigr)  =  0 $
  (or some other differential equality) in the periodicity cell $Y$, considered as an open set,
  without satisfying the same property in the entire space.
  Figure \ref{fig:bumps} presents a somewhat similar situation in one variable only.
  It shows a $ \mathbb Z $-periodic function $u$ such that $ u'' = -2 $ in $ Y = {]}0, 1{[} $
  but $ u'' \neq -2 $ in $ \mathbb R $;
  $ u''+2 $ is a sum of Dirac measures supported at $ \mathbb Z $.
\end{remark}
  
\begin{remark}
  The solution of problem (\ref{eq:var-form-displ-macrostrain}) is not unique;
  it contains an abitrary additive constant vector.
We need to tighten the space in order to obtain a unique solution.
We can do this by imposing a zero-average condition$\,$: $ u_A \in LP_0 \cap LP(A) $.
In this case, problem (\ref{eq:var-form-displ-macrostrain}) becomes
\begin{equation} \label{eq:var-form-displ-macrostr}
\left\{
\begin{matrix} \mbox{find } u_A \in LP_0 \cap LP(A) \mbox{ such that } \\
  \displaystyle\int_{Y} \bigl\langle \tens C\, \tens e(u_A) ,\tens e(v)\bigr\rangle\, dy =0
  \ \ \forall v \in H^1_{\#0} \,. 
\end{matrix} \right. 
\end{equation}

A roughly equivalent approach is to use a quotient space$\,$: $ LP(A) $ divided by the subspace
of constant vectors.
\end{remark}

Note that the classical Lax Milgram lemma does not apply to the above formulation
since the space $ LP(A) $ where the solution is to be found is a translation
of the space $ H^1_\# $ of test functions.
An adapted version of the Lax Milgram lemma was proven in \cite{Toader2011};
for the sake of self containedness we also state it as Lemma \ref{lemma:Lax-Milgram} below.
It ensures the existence and the uniqueness of the solution of the above formulation by taking
$ V = LP_0 $, $ V_0 = H^1_{\#0} $ and $ K = LP_0 \cap LP(A) $.

\begin{lemma}
\label{lemma:Lax-Milgram}
Let $V$ be a fixed Hilbert space, let $V_0$ be a closed subspace of $V$ and let $K = \gamma +V_0$
be a translation of $V_0$ (a closed affine manifold in $V$ ) where $\gamma$ is a fixed element of $V$.
Consider $a : V \times V \mapsto \RR$ a bilinear symmetric continuous form on $V$
which is coercive only on $V_0$, and consider $l : V \mapsto \RR$ a linear continuous form on V . Then
the problem
$$
\left\{
\begin{matrix} \mbox{find } u \in K \mbox{ such that } \\
a(u, v) = l(v)\,,\ \forall v \in V_0,
\end{matrix} \right. 
$$
has a unique solution.
\end{lemma}

Another variational formulation of the cellular problem is

\begin{equation}  \label{eq:var-form-displ-macrostrain-rep}
\left\{
\begin{matrix} \mbox{find } \ffi_A \in H^1_\# \mbox{ such that } \\
  \displaystyle\int_Y \langle \tens C\, \tens e(\ffi_A), \tens e(v)\rangle\, dy =
  - \bigl\langle A , \int_Y \tens C\, \tens e(v)\, dy \bigr\rangle\,,\ \ \forall v \in H^1_\# \,. 
\end{matrix} \right. 
\end{equation}

If we add a zero-average condition on $ \ffi_A $, the solution is unique\,:

\begin{equation}  \label{eq:var-form-displ-macrostrain-rep0}
\left\{
\begin{matrix} \mbox{find } \ffi_A \in H^1_{\#0} \mbox{ such that } \\
  \displaystyle\int_Y \langle \tens C\, \tens e(\ffi_A), \tens e(v)\rangle\, dy =
  - \bigl\langle A , \int_Y \tens C\, \tens e(v)\, dy \bigr\rangle\,,\ \ \forall v \in H^1_{\#0} \,. 
\end{matrix} \right. 
\end{equation}

The classical Lax Milgram Lemma ensures the existence of a solution $ \ffi_A $.

Formulation (\ref{eq:var-form-displ-macrostrain-rep0}) is frequently encountered in the literature,
often stated on the periodicity cell $Y$ rather than on the entire space $ \RR^3 $.
See e.g.\ \cite[Chapter I, section 2.2]{Bensoussan}, \cite[Chapter I, section 6.1]{Oleinik},
\cite[Definition 2.1]{Allaire1992}.

A variational formulation in strain arises naturally
(recall that the strain is the symmetric part of the gradient of the displacement)\,:

\begin{proposition} \label{thrm:var-form-strain-1}
The variational formulation in strain of the cellular problem
\begin{equation}  \label{eq:var-form-strain-1}
\left\{
\begin{matrix} \mbox{find } \tens e \in \TT^\perp \mbox{ such that } \\
  \displaystyle\media_Y \langle \tens C \tens e, \tens\mu \rangle\, dy =
  -\bigl\langle A , \media_Y \tens C \tens\mu\, dy \bigr\rangle\,,\ \ \forall \tens\mu \in \TT^\perp . 
\end{matrix} \right. 
\end{equation}
is equivalent to (\ref{eq:var-form-displ-macrostrain-rep0}).
\end{proposition}

\begin{proof}
Follows from the isomorphism between $\TT^\perp$ and $H^1_{\#0}$ (see Theorem \ref{thrm:isomorphism}).
\end{proof}

The proof of Proposition \ref{prop:equiv-uA-phiA} below is a mere exercise.

\begin{proposition}\label{prop:equiv-uA-phiA}
Formulations (\ref{eq:var-form-displ-macrostrain}) and
(\ref{eq:var-form-displ-macrostrain-rep}) are equivalent, the solutions $ u_A $ and $ \ffi_A $
being related by $ u_A(y) = Ay + \ffi_A(y) $.
Also, formulations (\ref{eq:var-form-displ-macrostr}) and (\ref{eq:var-form-displ-macrostrain-rep0})
are equivalent.
\end{proposition}

The variational formulation of the cellular problem
when the macroscopic stress $S $ is given in $ \M^s_3(\RR) $ is the following

\begin{equation} \label{eq:var-form-displ-macrostress-mean}
\left\{
\begin{matrix} \mbox{find } w_S \in LP_0 \mbox{ such that } \\
\displaystyle\media_Y \tens C\, \tens e(w_S)\, dy = S \\
\dint_{Y} \bigl\langle\tens C\, \tens e(w_S), \tens e(v)\bigr\rangle\, dy = 0 ,  \ \ 
\forall v \in H^1_{\#0}.
\end{matrix}
\right.
\end{equation}

Once again, one can apply
Lemma \ref{lemma:Lax-Milgram}, the adapted version of the Lax Milgram Lemma,
in order to ensure the existence and the uniqueness of the solution $w_S$;
here, $ V = LP_0 $, $ V_0 = H^1_{\#0} $ and $K$ is the set of functions in $ LP_0 $ satisfying
$ \displaystyle\media_Y \tens C\, \tens e(w_S)\, dy = S $.

An equivalent approach is to implement the average condition in the linear form$\,$:

\begin{equation}  \label{eq:var-form-displ-macrostress}
\left\{
\begin{matrix} \mbox{find } w_S \in LP_0 \mbox{ such that } \\
  \displaystyle\media_Y \langle\tens C\, \tens e(w_S) ,\tens e(v)\rangle\, dy =
  \bigl\langle S , \media_Y \tens e(v)\, dy \bigr\rangle,  \ \ 
\forall v \in LP_0 \,.
\end{matrix}
\right.
\end{equation}

\begin{proposition}\label{thrm:var-form-strain-0}
  The variational formulations (\ref{eq:var-form-displ-macrostress-mean}) and
  (\ref{eq:var-form-displ-macrostress}) are equivalent.
\end{proposition}

\begin{proof}
  Take $ w_S $ satisfying (\ref{eq:var-form-displ-macrostress-mean}) and let $ v\in LP_0 $;
  we want to prove the equality in (\ref{eq:var-form-displ-macrostress}).
  Any $ v\in LP_0 $ can be written as a sum $ v(y) = By + \bar v(y) -
  \displaystyle \media_Y By\,dy $ for some symmetric matrix $B$ and for some $ \bar v \in H^1_{\#0} $.
  Then
  $$ \media_Y \bigl\langle\tens C\, \tens e(w_S), \tens e(v)\bigr\rangle\, dy =
  \media_Y \bigl\langle\tens C\, \tens e(w_S), B \bigr\rangle\, dy +
  \media_Y \bigl\langle\tens C\, \tens e(w_S) ,\tens e(\bar v) \bigr\rangle\, dy =
  \langle B, S \rangle   $$
  which proves $ w_S $ satisfies (\ref{eq:var-form-displ-macrostress}).

  Conversely, take $ w_S $ satisfying (\ref{eq:var-form-displ-macrostress}).
  We begin by choosing
  the test function $ v(y) = By - \displaystyle \media_Y By\,dy $
  for an arbtrarily chosen symmetric matrix $B$.
  Then the equality in (\ref{eq:var-form-displ-macrostress}) implies
  $ \displaystyle \media_Y \tens C\, \tens e(w_S)\, dy = S $.
  Now choose $ v\in H^1_{\#0} $ and, since the gradient of any
  periodic function has zero average, we obtain the third line in
  (\ref{eq:var-form-displ-macrostress-mean}).
\end{proof}

Again, a variational formulation in strain arises naturally:

\begin{proposition}\label{thrm:var-form-strain-2}
The variational formulation in strain of the cellular problem
\begin{equation}  \label{eq:var-form-strain-2}
\left\{
\begin{matrix} \mbox{find } \tens e \in \TT(0)^\perp  \mbox{ such that } \\
  \displaystyle\media_Y \langle \tens C\, \tens e, \tens\mu \rangle\, dy =
  \bigl\langle S , \media_Y  \tens\mu\, dy \bigr\rangle\,,\ \ \forall \tens\mu \in \TT(0)^\perp . 
\end{matrix} \right. 
\end{equation}
is equivalent to (\ref{eq:var-form-displ-macrostress}).
\end{proposition}

\begin{proof}
Follows from the isomorphism between $\TT(0)^\perp$ and $LP_0$ (see Theorem \ref{thrm:LP-isomorphic-TT_0-perp}).
\end{proof}

\begin{proposition}\label{prop:var-form-3}
  The variational formulations (\ref{eq:var-form-displ-macrostr}) and
  (\ref{eq:var-form-displ-macrostress-mean}) are equivalent, in the following sense\,:
  Consider a macroscopic strain $A$; suppose $ u_A $ satisfies formulation
  (\ref{eq:var-form-displ-macrostr});
  define $ S = \displaystyle \media_Y \tens C\, \tens e(u_A) $ (or, equivalently, $ S = \tens C^H\! A $);
  then $ u_A $ is a solution of (\ref{eq:var-form-displ-macrostress-mean}).
  Conversely, let $S$ be a macroscopic stress; suppose $ w_S $ satisfies formulation
  (\ref{eq:var-form-displ-macrostress-mean});
  define $ A= \displaystyle \media_Y \tens e(w_S) $ (or, equivalently, $ A = \tens D^H\! S $);
  then $ w_S $ is a solution of problem (\ref{eq:var-form-displ-macrostr}).
\end{proposition}

\begin{proof}
  Both implications stated above are rather obvious due to the choice of the matrices $A$ and $S$.
  The second implication relies also on Remark \ref{rem:LPA-avrg}.
\end{proof}

Problem (\ref{eq:var-form-displ-macrostr}) is still equivalent with the following formulation
in stress, provided $ S = \tens C^H A $.

\begin{equation}\label{eq:var-form-macrostress}
\left\{ \begin{matrix} \mbox{find } \tens\sigma\!{}_S \in \TT(S) \mbox{ such that } \\
\dint_{Y} \bigl\langle \tens D \tens\sigma\!{}_S , \tens\mu \bigr\rangle\, dy = 0 \ \ 
\forall \tens\mu  \in \TT(0) ,
\end{matrix}
\right.
\end{equation}
where the existence and the uniqueness of the solution $\tens\sigma\!{}_S$ is guaranteed by
the adapted version of the Lax-Milgram Lemma \ref{lemma:Lax-Milgram} above,
applied with $ V=\TT $, $ V_0= \TT(0) $ and $ K=\TT(S) $.

\begin{proposition}\label{prop:var-form-4}
  The variational formulations (\ref{eq:var-form-displ-macrostress-mean}) and
  (\ref{eq:var-form-macrostress}) are equivalent.
\end{proposition}

\begin{proof}
  Take $ w_S $ satisfying (\ref{eq:var-form-displ-macrostress-mean}) and define
  $ \tens\sigma\!{}_S = \tens C \tens e(w_S) $.
  Then $ \tens\sigma\!{}_S $ is periodic and its average is $S$.
  Also, $ \div\tens\sigma\!{}_S =0 $ since the formulations
  (\ref{eq:cell-pb}), (\ref{eq:var-form-displ-macrostrain}) versus (\ref{eq:var-form-displ-macrostr})
  and (\ref{eq:var-form-displ-macrostress-mean}) are equivalent up to an additive constant vector
  (see Propositions \ref{prop:var-form-1} and \ref{prop:var-form-3} above).
  Therefore, $ \tens\sigma\!{}_S $ belongs to $ \TT(S) $.
  On the other hand, recalling that $ \tens D $ is the inverse of the rigidity tensor $ \tens C $,
  one has $ \displaystyle \media_{Y} \bigl\langle \tens D \tens\sigma\!{}_S ,\tens\mu \bigr\rangle\, dy =
  \media_{Y} \bigl\langle \tens e(w_S) , \tens\mu \bigr\rangle\, dy $.
  Due to Theorem \ref{thrm:prod-avrg}, the above quantity is equal to $ \displaystyle
  \bigl\langle \media_Y \tens e(w_S) \,dy, \media_Y \tens\mu\, dy \bigr\rangle $ which is zero
  because $ \tens\mu\in\TT(0) $.

  Conversely, let $ \tens\sigma\!{}_S $ be the solution of (\ref{eq:var-form-macrostress}).
  Due to Theorem \ref{thrm:Donati.2}, there exists $ w_S \in LP $ such that $ \tens D \tens\sigma\!{}_S
  = \tens e(w_S) $ thus $ \tens\sigma\!{}_S = \tens C \tens e(w_S) $.
  We can use the arbitrary additive constant vector mentioned in Theorem \ref{thrm:Donati.2}
  and choose a function $ w_S $ having zero average.
  On the other hand, since $ \tens C \tens e(w_S) \in\TT_S \subset \TT $ and
  taking a test function $ v\in H^1_{\#0} \subset LP_0 \subset LP $, Theorem \ref{thrm:prod-avrg} implies
  $ \displaystyle \media_Y \bigl\langle \tens C \tens e(w_S), \tens e(v) \bigr\rangle =
  \bigl\langle \media_Y \tens C \tens e(w_S) \,dy, \media_Y \tens e(v)\,dy \bigr\rangle $.
  Again, because the gradient of any periodic function has zero average,
  $ w_S $ satisfies (\ref{eq:var-form-displ-macrostress-mean}).
\end{proof}

\begin{figure}[ht] \centering
  \psfrag{18}{(\ref{eq:var-form-displ-macrostrain})}
  \psfrag{19}{(\ref{eq:var-form-displ-macrostr})}
  \psfrag{20}{(\ref{eq:var-form-displ-macrostrain-rep})}
  \psfrag{21}{(\ref{eq:var-form-displ-macrostrain-rep0})}
  \psfrag{22}{(\ref{eq:var-form-strain-1})}
  \psfrag{23}{(\ref{eq:var-form-displ-macrostress-mean})}
  \psfrag{24}{(\ref{eq:var-form-displ-macrostress})}
  \psfrag{25}{(\ref{eq:var-form-strain-2})}
  \psfrag{26}{(\ref{eq:var-form-macrostress})}
  \psfrag{prop15}{Prop.\ \ref{thrm:var-form-strain-1}}
  \psfrag{prop16}{Prop.\ \ref{prop:equiv-uA-phiA}}
  \psfrag{prop17}{Prop.\ \ref{thrm:var-form-strain-0}}
  \psfrag{prop18}{Prop.\ \ref{thrm:var-form-strain-2}}
  \psfrag{prop19}{Prop.\ \ref{prop:var-form-3}}
  \psfrag{prop20}{Prop.\ \ref{prop:var-form-4}}
  \psfrag{u}{unknown displacement}
  \psfrag{f}{unknown $\ffi$}
  \psfrag{e}{unknown strain}
  \psfrag{s}{unknown stress}
  \psfrag{A}{datum $A$}
  \psfrag{S}{datum $S$}
  \includegraphics[width=119mm]{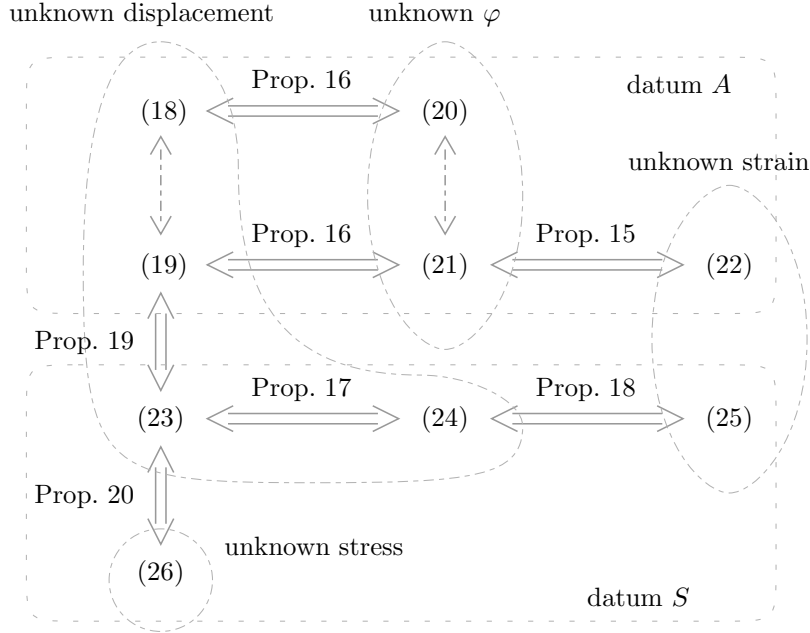}
  \caption{Relationships between different variational formulations of the cellular problem}
  \label{fig:implications}
\end{figure}

Figure \ref{fig:implications} shows a diagram of the equivalences between the variational formulations
presented in the current section.
Problems (\ref{eq:var-form-displ-macrostrain}) and (\ref{eq:var-form-displ-macrostrain-rep})
have non-unique solutions; by adding a zero-average condition we obtain problems
(\ref{eq:var-form-displ-macrostr}) and (\ref{eq:var-form-displ-macrostrain-rep0}) with unique solution.
Joint to each equivalence sign, the proposition stating it is specified.
Formulations (\ref{eq:var-form-displ-macrostrain}), (\ref{eq:var-form-displ-macrostr}), 
(\ref{eq:var-form-displ-macrostrain-rep}), (\ref{eq:var-form-displ-macrostrain-rep0})
and (\ref{eq:var-form-strain-1}) have the macroscopic strain $A$ as datum.
Formulations (\ref{eq:var-form-displ-macrostress-mean}), (\ref{eq:var-form-displ-macrostress}),
(\ref{eq:var-form-strain-2}) and (\ref{eq:var-form-macrostress}) have the macroscopic stress $S$ as datum.
In formulations (\ref{eq:var-form-displ-macrostrain}), (\ref{eq:var-form-displ-macrostr}), 
(\ref{eq:var-form-displ-macrostress-mean}) and (\ref{eq:var-form-displ-macrostress})
the unknown is the microscopic displacement.
In formulations (\ref{eq:var-form-displ-macrostrain-rep}) and (\ref{eq:var-form-displ-macrostrain-rep0})
the uknown is $ \ffi $ which is the periodic part of the microscopic displacement.
In formulations (\ref{eq:var-form-strain-1}) and (\ref{eq:var-form-strain-2}) the unknown is the
microscopic strain.
Finally, in formulation (\ref{eq:var-form-macrostress}) the unknown is the microscopic stress.

\section{The cellular problem formulated as minimization problem}
\label{sec:cell-pb-min}

In the sequel we extend to the periodic framework the results in \cite[sections 4 and 5]{Ciarlet-G-K2011},
in order to obtain formulations of the cellular problem (\ref{eq:cell-pb})
in displacement, stress and strain.
We relate these formulations to the variational formulations (\ref{eq:var-form-displ-macrostrain}),
(\ref{eq:var-form-displ-macrostrain-rep}), (\ref{eq:var-form-displ-macrostress-mean}),
(\ref{eq:var-form-displ-macrostress}) and (\ref{eq:var-form-macrostress}) in section \ref{sec:var-form}.
Recall that these variational formulations are all equivalent in the sense that $ u_A = w_S $
and $ \tens\sigma\!{}_S = \tens C\, \tens e(w_S) $, provided $ S = \tens C^H A $.

\begin{theorem}[formulation in displacement]
\label{thrm:displ}
Given  a macroscopic strain $A$, there exists a unique vector field $ u_A $ in $ LP_0 \cap LP(A) $
that satisfies
$$
\J (u_A) = \displaystyle{\inf_{u \in LP_0\cap LP(A)} } \J(u)
$$
where  
$$
\J(u) = \dfrac{1}{2} \dint_Y \langle \tens C \tens e(u), \tens e(u) \rangle\, dy \hbox{ for all } u\in LP_0\cap LP(A) .
$$
Moreover, $u_A$ is the solution of the variational formulation (\ref{eq:var-form-displ-macrostr}).
\end{theorem}

\begin{proof}
Given $ u \in LP(A) $ the Gâteaux derivative of $ \J $ is
$
D\J(u)(v) = \dint_Y \langle \tens C\, \tens e(u), \tens e(v) \rangle\, dy
$
for all $v \in H^1_{\#0} $. 
The condition $D \J(u)(v) =0$ $\forall v \in H^1_{\#0} $ is equivalent to the variational
formulation  (\ref{eq:var-form-displ-macrostr}).
Consequently there exists a unique $ u_A \in LP_0 \cap LP(A) $ with
the property that $ D \J(u_A)(v)=0 $. The functional $ \J $ attains its infimum in $ u_A $.
\end{proof}

Note that the value of $\J$ does not change if one adds a constant vector to $ u_A $, so it easy to
extend $\J$ to $ LP(A) $.

Another formulation in displacement is needed in order to minimize the function over a Hilbert space
$ H^1_{\#0} $ rather than over an afine space $ LP_0 \cap LP(A) $.
The next theorem gives therefore a formulation in displacement, equivalent to the above one.  

\begin{theorem}[equivalent formulation in displacement]
\label{thrm:displ.2}
1) Given a macroscopic strain $A$, there exists a unique vector field $\ffi_A$ in $ H^1_{\#0} $
satisfying
$$
J (\ffi_A) = \displaystyle{\inf_{\ffi \in H^1_{\#0}} } J(\ffi)
$$
where  
$$
J(\ffi) = \dfrac{1}{2} \dint_Y \langle \tens C \tens e(\ffi), \tens e(\ffi) \rangle\, dy
+ \langle A, \dint_Y \tens C \tens e(\ffi)\, dy \rangle \hbox{ for all } \ffi\in H^1_{\#0} .
$$
Moreover, $\ffi_A$ is the solution of the variational formulation
(\ref{eq:var-form-displ-macrostrain-rep0}).

2) Given a macroscopic stress $S$, there exists a unique vector field $w_S$ in $LP_0$ that satisfies
$$
K (w_S) = \displaystyle{\inf_{u \in LP_0} } K(u)
$$
where  
$$
K(u) = \dfrac{1}{2} \dint_Y \langle \tens C \nabla_s u ,\nabla_s u \rangle\, dy
- \dint_Y \langle S ,  \nabla_s u \rangle \, dy \hbox{ for all } u\in LP_0 \,.
$$
Moreover, $w_S$ is the solution of the variational formulation (\ref{eq:var-form-displ-macrostress}).

\end{theorem}

\begin{proof}
1) Given $ \ffi \in H^1_{\#0} $, the G\^ateaux derivative of $J$ is
$
D\,J(\ffi)(\psi) = \dint_Y \langle \tens C\, \tens e(\ffi), \tens e(\psi) \rangle\, dy
+ \langle A, \dint_Y \tens C \tens e(\psi)\, dy \rangle
$
for all $ \psi \in H^1_{\#0} $.  
The derivative $D J(\ffi)(\psi) =0 $ $\forall \psi \in H^1_{\#0} $ is equivalent to the variational
formulation  (\ref{eq:var-form-displ-macrostrain-rep0}).
Consequently there exists a unique $ \ffi \in H^1_{\#0} $ with the property that $ D \J(\ffi)(\psi)=0 $.
Therefore the functional $ J $ attains its infimum in $ \ffi_A $.

2) The proof is analogous with the one in 1)\,: the G\^ateaux derivative of $K$ being zero
is equivalent to the variational formulation (\ref{eq:var-form-displ-macrostress}).
\end{proof}

\begin{theorem}[formulation in stress] \label{thrm:form-stress}
Given a macroscopic stress $ S \in \M_3^s(\RR) $,
there exists a unique tensor field $\tens\sigma\!{}_S \in \TT(S)$ that satisfies
\begin{equation} \label{eq:func-g}
g(\tens\sigma\!{}_S) = \displaystyle{\inf_{\tens\sigma \in \TT(S)} g(\tens\sigma)}
\end{equation}
\noindent where
$$
g(\tens\sigma) = \dfrac{1}{2} \dint_Y \langle \tens D \tens\sigma, \tens\sigma \rangle\, dy \hbox{ for all } \tens\sigma \in \LL^2_{\# s} 
$$
and $ \tens D = \tens C^{-1} $ is the compliance tensor.

\noindent  Moreover, $ \tens\sigma\!{}_S $ is the solution of the variational formulation
(\ref{eq:var-form-macrostress}) and $ \tens D \tens\sigma\!{}_S = \nabla\!{}_s u_A $,
where $ u_A $ is the solution of (\ref{eq:var-form-displ-macrostr}).
\end{theorem}

\begin{proof}
For $\tens\sigma \in \TT(S)$, the Gâteaux derivative of $g$ is
$
D g(\tens\sigma)(\tens\mu) = \dint_Y \langle \tens D \tens\sigma , \tens\mu \rangle\, dy 
$
for all $\tens\mu \in \TT(0)$.
Then the derivative $D g(\tens\sigma)(\tens\mu)=0 \,  \forall \tens\mu  \in  \TT(0)$ is equivalent to
the above variational formulation in stress of the cellular problem (\ref{eq:var-form-macrostress})
and therefore the solution exists and is unique so $\tens\sigma = \tens\sigma\!{}_S$.

On the other hand, if $\dint_Y \langle \tens D \tens\sigma\!{}_S , \tens\mu \rangle\, dy = 0
\quad \forall \tens\mu  \in \TT(0) $, applying the
extension of the Donati's Theorem \ref{thrm:Donati.2} it turns out that $\tens D \tens\sigma\!{}_S = \nabla\!{}_s v$ for a function $v \in LP$ 
which is precisely $v=u_A$ for $ A= \tens D^H S $;
recall that $ \tens D^H = (\tens C^H)^{-1} $ is the homogenized compliance tensor.
\end{proof}

Theorem \ref{thrm:form-strain} below states the cellular problem as a minimization problem on
the space $ \TT(0)^\perp $, defined in formula (\ref{eq:def-TT-0-ort}).
It is inspired in \cite[Theorem 4.3]{Ciarlet-G-K2011} where the minimization is performed on
a similar space $ \MM^\perp $.
The definition of the space $ \MM $, stated in \cite[Theorem 3.2]{Ciarlet-G-K2011}, is based on
a space $V$ which contains a Dirichlet boundary condition.
Consequently, $ \MM $ contains an implicit Neumann condition which casts a Dirichlet-like
condition on elements of $ \MM^\perp $.
Similarly, the zero-average condition appearing in the definition (\ref{eq:def-TT-S}) of
the space $ \TT(0) $ can be viewed as a homogeneous Neuman condition.
This casts a condition on elements of $ \TT(0)^\perp $ which has the nature of a Dirichlet
boundary condition.
See also \cite[Remark 2]{BT2010}, \cite[Remark 1]{Toader2011} and \cite[end of Section 3]{BL2011}.

\begin{theorem}[formulation in strain] \label{thrm:form-strain}
Each of the minimization problems below has a unique solution (a microscopic strain field)\,:  

\noindent 1) For a given macroscopic strain $A$, find $\bar {\tens e} \in \TT^\perp $ such that
$$
j(\bar {\tens e}) = \displaystyle{\inf_{\tens e \in \TT^\perp} j(\tens e)},
$$
where $j : \LL^2_{\# s} \mapsto \RR$ is defined by
$
\displaystyle{ j(\tens e) = \dfrac{1}{2} \dint_Y \langle \tens C \tens e, \tens e \rangle\, dy
  + \langle A, \dint_Y\tens C \tens e\, dy \rangle} ;
$

\noindent 2) For a given macroscopic stress $S$, find $\tilde {\tens e} \in \TT(0)^\perp $ such that
$$
k(\tilde {\tens e}) = \displaystyle{\inf_{\tens e \in \TT(0)^\perp} k(\tens e)},
$$
where $k : \LL^2_{\# s} \mapsto \RR$ is defined by
$
\displaystyle{k(\tens e) = \dfrac{1}{2} \dint_Y \langle \tens C \tens e, \tens e \rangle\, dy
  - \langle S, \dint_Y \tens e\, dy \rangle .}
$

\noindent Moreover, the above two minimization problems are equivalent, in the following sense.
Consider a macroscopic strain $A$; take the solution $ \bar {\tens e} $ of problem 1) above
(associated to $A$); define $ S = \displaystyle{\media_Y \tens C\bar {\tens e}} $ (or, equivalently,
$ S = \tens C^H\! A $);
take the solution $ \tilde {\tens e} $ of problem 2) above (associated to $S$);
then $ \bar {\tens e} = \tilde {\tens e} -A $ and $ j ( \bar {\tens e} ) = k ( \tilde {\tens e} ) $.
Conversely, let $S$ be a macroscopic stress; take the solution $ \tilde {\tens e} $ of problem 2) above
(associated to $S$); define $ A= \displaystyle{\media_Y \tilde {\tens e}} $ (or, equivalently,
$ A = \tens D^H\! S $);
take the solution $ \bar {\tens e} $ of problem 1) above (associated to $A$);
then $ \tilde {\tens e} = \bar {\tens e} +A $ and $ k ( \tilde {\tens e} ) = j ( \bar {\tens e} ) $.
\end{theorem}

\begin{proof}
1) Deriving the function $j$ one obtains 
$
D j(\tens e)(\tens\mu) = \dint_Y \langle \tens C \tens e , \tens\mu \rangle\, dy
+ \langle A, \dint_Y\tens C \tens\mu\, dy \rangle
$
for all $ \tens\mu \in \TT^\perp $.
Then $ D j(\tens e)(\tens\mu) =0 $ for all $ \tens\mu \in \TT^\perp $, is equivalent to the formulation (\ref{eq:var-form-strain-1}).
From Proposition \ref{thrm:var-form-strain-1}  it follows that the solution $\bar {\tens e}$ is unique and it verifies 
$\bar{\tens e}= \nabla\!{}_s \ffi_A$,
thus $\bar {\tens e} = \nabla\!{}_s u_A - A$.

2) The G\^ateaux derivative of $k$ is  
$
D k(\tens e)(\tens\mu) = \dint_Y \langle \tens C \tens e , \tens\mu \rangle\, dy
- \langle S, \dint_Y \tens\mu\, dy \rangle
$
for all $\tens\mu \in \TT(0)^\perp $, and $D k(\tens e)(\tens\mu) =0 $ for all $\tens\mu \in \TT(0)^\perp $,
is equivalent to the formulation (\ref{eq:var-form-strain-2}).

Thus, Proposition \ref{thrm:var-form-strain-2} implies that the solution $ \tilde {\tens e} $
is unique and verifies $ \tilde {\tens e} = \nabla\!{}_s u_A$.

The integral conditions follow from (\ref{eq:A}) and (\ref{eq:S}).
\end{proof}

The state of a body is characterized by the internal stress and the local displacement.
Its energy is the product between the stress tensor and the strain tensor,
defined as the symmetric gradient of the displacement.
We introduce an energy operator $ \Lambda : \LL^2_{\# s} \mapsto (LP)^\star $, where $ (LP)^\star $
is the dual of the Hilbert space $LP$, defined by 
$$
\langle \Lambda \tens\sigma , v \rangle = \int_Y \langle \tens\sigma, \nabla\!{}_s v \rangle\, dy , \forall v \in LP.
$$
Consider $ S \in \M_3^s(\RR) $ a stress matrix.
Consider the particular energy functional defined as
\begin{equation} \label{eq:def-eta-S}
  \eta_S(v) = \langle \Lambda S, v \rangle = \int_Y \langle S,  \nabla\!{}_s v \rangle\, dy \,.
\end{equation}
Define the indicator function $ \I_S: (LP)^\star \mapsto \RR \cup \{ +\infty \}$ defined by
\begin{equation} \label{eq:def-IS}
\I_S(\eta) = \displaystyle{\left\{\begin{matrix}
0 \mbox{ if } \eta =\eta_S\,, \\
+\infty \mbox{ if } \eta \neq \eta_S \,.
  \end{matrix} \right. }
\end{equation}

Note that $ \I_S(\Lambda\tens\sigma) = \displaystyle{\left\{\begin{matrix}
  0 \mbox{ if } \langle \Lambda\tens\sigma, v\rangle = \langle \Lambda S, v\rangle\,, \\
  +\infty \mbox{ if } \langle \Lambda\tens\sigma, v\rangle \neq \langle \Lambda S, v\rangle\,.
  \end{matrix}\right.} $
$ \I_S(\Lambda\tens\sigma) $ is zero if and only if
$ \displaystyle \int_Y \langle \tens\sigma, \nabla\!{}_s v \rangle\, dy =
\int_Y \langle S, \nabla\!{}_s v \rangle\, dy $, $ \forall v\in LP $.

Theorem \ref{thrm:minimization.in.stress} below transforms the stress formulation
in Theorem \ref{thrm:form-stress} in a minimization over all symmetric stresses
in $ \LL^2_{\# s} $ of a function that contains implicitly the variational formulation
(\ref{eq:var-form-macrostress}).
As we shall see during the proof, the second term in (\ref{eq:G}), $ \I_S(\Lambda\tens\sigma) $,
ensures that the stress $ \tens\sigma $ (minimizer of $G$) belongs to $ \TT(S) $.

\begin{theorem} \label{thrm:minimization.in.stress} 
The energy operator $\Lambda$ belongs to $ \L(\LL^2_{\# s}, (LP)^\star) $.
The functional $ \I_S $ above defined is proper, convex and lower semicontinuous and verifies
$ \I_S(\Lambda\tens\sigma)= I_{\TT(S)}(\tens\sigma) $ where $ I_{\TT(S)} $ denotes the indicator function
of the set $ \TT(S) $.
The function $g$ defined in (\ref{eq:func-g}) is proper, convex and lower semicontinous.

\noindent Defining the function $ G : \LL^2_{\# s} \mapsto \RR \cup \{ +\infty \}$ by 
\begin{equation} \label{eq:G}
G(\tens\sigma )= g(\tens\sigma ) + \I_S(\Lambda \tens\sigma ), \quad \forall \tens\sigma \in \LL^2_{\# s} ,
\end{equation}
the minimization problem in Theorem \ref{thrm:form-stress} is equivalent to the following one\,:
\begin{equation}
\label{eq:primal-problem-G}
 \displaystyle{\inf_{\tens\sigma \in \LL^2_{\# s}} G(\tens\sigma)},
\end{equation}
in the sense that the minimum value is the same and is attained for the same $\tens\sigma$\,: 
$ \displaystyle{\inf_{\tens\sigma \in \TT(S)} g(\tens\sigma)} \Leftrightarrow
\displaystyle \inf_{\tens\sigma \in \LL^2_{\# s}} G(\tens\sigma) $.
\end{theorem}

\begin{proof}
For any $ \tens\sigma \in \LL^2_{\# s} $ the functional $ \Lambda \tens\sigma : LP \mapsto \RR $ defined by
$ \displaystyle \langle \Lambda \tens\sigma, v \rangle =  \int_Y \langle \tens\sigma, \nabla\!{}_s v \rangle\, dy $
is linear and continuous, therefore $ \Lambda \tens\sigma $ belongs to $ (LP)^\star $ and the definition of
$\Lambda$ is consistent.
Moreover, for all $ \tens\sigma \in \LL^2_{\# s} $
 $$
 \| \Lambda\tens\sigma\|_{(LP)^\star} = \displaystyle{\sup_{v\in LP} \dfrac{\int_Y \langle \tens\sigma, \nabla\!{}_s v \rangle}{\| v\|_{H^1_\#}}} \le 
 \displaystyle{\sup_{v\in LP} \dfrac{\int_Y \langle \tens\sigma, \nabla\!{}_s v \rangle}{\|\nabla\!{}_s v\|_{\LL^2_{\# s}}}} = 
 \| \tens\sigma \|_{\LL^2_{\# s}},
 $$
and consequently $ \Lambda \in\L(\LL^2_{\# s}, (LP)^\star) $.

The functional $\I_S$ is the indicator function of $ \{ \eta_S \} \subset (LP)^\star $
and therefore it is proper,
convex and lower-semicontinuous (since $ \{ \eta_S \} $ is closed in $ (LP)^\star$ ). 
 
\noindent Let us prove that $ \I_S(\Lambda\tens\sigma)=  I_{\TT(S)}(\tens\sigma) $,
which is equivalent to $ \Lambda\tens\sigma = \eta_S $, that is,
$$
\langle \Lambda \tens\sigma, v \rangle = \langle S, \nabla\!{}_s v \rangle ,\,\forall v \in LP  \quad \Leftrightarrow
$$
$$
\langle \tens\sigma , \nabla\!{}_s v \rangle=\langle S, \nabla\!{}_s v \rangle ,\, \forall v \in LP  \quad \Leftrightarrow
$$
$$
\langle \tens\sigma - S, \nabla\!{}_s v \rangle = 0 ,\, \forall v \in LP .
$$
The above equality contains the mean condition on $\tens\sigma $ and the zero divergence condition
and is thus equivalent to $\tens\sigma \in \TT(S)$.
We conclude that $ \I_S(\Lambda\tens\sigma)=  I_{\TT(S)}(\tens\sigma)$,
for all $\tens\sigma \in \LL^2_{\# s} $.

The function $g$ is proper, convex (quadratic in $\tens\sigma$) and continuous for the norm
$ \| \cdot \|_{\LL^2_{\#}} $.
The minimization problem in Theorem \ref{thrm:form-stress} is equivalent to 
$$
\displaystyle{\inf_{\tens\sigma \in \LL^2_{\# s}} g(\tens\sigma)+ I_{\TT(S)}(\tens\sigma) }
$$
and the last assertion in the theorem is a consequence of
$ \I_S(\Lambda\tens\sigma) $ being equal to $ I_{\TT(S)}(\tens\sigma) $.
\end{proof}

Theorem \ref{thrm:L-F-transf} below
describes the Legendre-Fenchel transforms for the functions $g$ and $\I_S$
introduced in Theorem \ref{thrm:form-stress} and 
respectively in Theorem \ref{thrm:minimization.in.stress} above.
We refer the reader to \cite{ET1972} and \cite{Brezis2010} for the theory and properties
of the Legendre-Fenchel transform; we merely state its definition\,:

\begin{definition}\label{def:L-F-transf}
  Let $X$ be a normed vector space and let $ \theta : X \to \RR\cup\{+\infty\} $ be a proper function.
  The Legendre-Fenchel transform of $\theta$ is the function
  $ \theta^\star : X^\star \to \RR\cup\{+\infty\} $ defined by
  $$ \theta^\star(e) = \sup_{x\in X}
  \bigl\{\langle e,x\rangle - \theta(x) \bigr\}\,, \mbox{ with } e\in X^\star\,. $$
\end{definition}

\begin{theorem}  \label{thrm:L-F-transf}
 The functions $g$ defined in (\ref{eq:func-g}) and $ \I_S $ defined in (\ref{eq:def-IS})
 have their Legendre-Fenchel transforms $ g^\star : \LL^2_{\# s}  \mapsto \RR $ and
 $ \I_S^\star : LP \mapsto \RR $ respectively given by
$$
g^\star (e) = \dfrac{1}{2} \dint_Y \langle \tens C e, e \rangle\, dy\, , \ \forall e\in  \LL^2_{\# s} \,,
$$
$$
\I_S^\star (v) = \eta_S(v), \,\forall v\in LP\,, 
$$
\noindent where $ \eta_S $ is the functional defined in (\ref{eq:def-eta-S}).
\end{theorem}

The proof is straightforward from the Definition \ref{def:L-F-transf} of the Legendre-Fenchel transform.

\section{Dual problems for the formulation in stress}
\label{sec:dual prob}

We introduce two dual problems to the minimization problem in stress,
according to the duality theory, see \cite[Chapter 6]{ET1972}.
A dual problem in strain and a dual problem in displacement will be presented,
generalizing the results in \cite[Sections 6 and 7]{Ciarlet-G-K2011} to the periodicity context.
Two Lagrangeans are introduced\,: $\L (\tens\sigma, \tens e) =\int_Y \langle \tens e , \tens\sigma \rangle\, dy - g^\star (\tens e) + \I_S(\Lambda \tens\sigma)$, for 
the stress-strain duality, and $\tilde \L (\tens\sigma, v) = \dfrac 1 2 \int_Y \langle \tens D \tens\sigma , \tens\sigma  \rangle\, dy +
 \langle \Lambda \tens\sigma , v  \rangle - \I_S^\star (v)$, for the stress-displacement duality.
 
We begin by presenting in detail the stress-strain Lagrangean $\L (\tens\sigma, \tens e)$ (Theorem \ref{thrm:Lagrangean-stress-strain} and 
Lemma \ref{thrm:dual-problem}) and afterwards the stress-displacement Lagrangean 
$\tilde \L (\tens\sigma, v)$ (Theorem \ref{thrm:Lagrangean-stress-displ} and 
Lemma \ref{thrm:dual-problem-displ}).

The minimization problem in stress (\ref{eq:primal-problem-G}) (primal problem) has a dual problem in strain, as we shall prove in the sequel.
A Lagrangian associated to the dual problem is introduced in the following theorem and its saddle-point corresponds to the solutions 
of the primal and dual problems, respectively (\ref{eq:primal-problem-G}) and (\ref{eq:dual_problem}).

\begin{theorem}
\label{thrm:Lagrangean-stress-strain}
Consider the Lagrangian $\L : \LL^2_{\# s} \times \LL^2_{\# s} \mapsto \RR \cup \{ + \infty \}$ defined by
\begin{equation}
\label{eq:Lagrangean-stress-strain}
\L (\tens\sigma, \tens e) = \int_Y \langle \tens e , \tens\sigma \rangle\, dy - g^\star (\tens e) + 
\I_S(\Lambda \tens\sigma) \quad \forall (\tens\sigma , \tens e) 
\in \LL^2_{\# s} \times \LL^2_{\# s}
\end{equation}
Then
$$
\displaystyle{\inf_{\tens\sigma \in \LL^2_{\# s}} G(\tens\sigma)=\inf_{\tens\sigma \in \LL^2_{\# s}}\sup_{\tens e \in \LL^2_{\# s}} \L (\tens\sigma, \tens e) = 
\L (\tens\sigma\!{}_S, \tens e(u_A)) =
\sup_{\tens e \in \LL^2_{\# s}}\inf_{\tens\sigma \in \LL^2_{\# s}} \L (\tens\sigma, \tens e)},
$$
where $\tens\sigma\!{}_S$ is the solution of both minimization problems in Theorem \ref{thrm:form-stress} and Theorem 
\ref{thrm:minimization.in.stress}  
and $\tens e (u_A)= \nabla_s u_A$ is the solution 
in Theorem \ref{thrm:form-strain}, provided $ S = \tens C^H A $. 
\end{theorem}

\begin{proof}
  The Fenchel-Moreau Theorem, presented in \cite{Fenchel1949} and \cite{Moreau1970},
  states that $g^{\star\star} = g$,
  where $g^\star$ is the Legendre-Fenchel transform of $g$.
Using the definition of Legendre-Fenchel of $g^\star$, $G$ may be written as follows\,:
$$
G(\tens\sigma)= g^{\star\star} (\tens\sigma) + \I_S(\Lambda \tens\sigma) = \sup_{\tens e \in \LL^2_{\# s}}
 \{ \int_Y \langle \tens e , \tens\sigma \rangle -g^\star (\tens e) \} + \I_S (\Lambda \tens\sigma) = \sup_{\tens e \in \LL^2_{\# s}} 
 \L (\tens\sigma, \tens e),
$$
for every $\tens\sigma \in \LL^2_{\# s}$.
Defining $G^\star : \LL^2_{\# s} \mapsto \RR\cup \{-\infty \}$ by
\begin{equation}
\label{eq:G-star}
G^\star (\tens e) = \inf_{\tens\sigma \in \LL^2_{\# s}} \bigl\{ \int_Y \langle \tens e , \tens\sigma \rangle\, dy
+ \I_S (\Lambda \tens\sigma) \bigr\} - g^\star(\tens e) ,
\end{equation}
the dual problem of (\ref{eq:primal-problem-G}) is
\begin{equation}
\label{eq:dual_problem}
 \displaystyle{\sup_{\tens e \in \LL^2_{\# s}} G^\star (\tens e)}.
\end{equation}
Lemma \ref{thrm:dual-problem} below ensures that
$$
G(\tens\sigma\!{}_S) = \inf_{\tens\sigma \in \LL^2_{\# s}} G(\tens\sigma) = \sup_{\tens e \in \LL^2_{\# s}} G^\star (\tens e) = 
G^\star (\tens e(u_A))
$$
and since one can verify that $G(\tens\sigma\!{}_S) = \L (\tens\sigma\!{}_S, \tens e (u_A))$, the conclusion of the theorem follows.
In order to confirm the above equality, recall that $\tens\sigma\!{}_S$ belongs to $\TT(S)$ and according to Theorem \ref{thrm:minimization.in.stress}, $h(\Lambda \tens\sigma\!{}_S) =0$. 
Using Theorem \ref{thrm:L-F-transf} one obtains
$$
\L (\tens\sigma\!{}_S, \tens e (u_A))= \int_Y \langle \tens e(u_A) , \tens\sigma\!{}_S \rangle\, dy - g^\star (\tens e(u_A)) + \I_S(\Lambda \tens\sigma\!{}_S) =
\int_Y \langle \tens e(u_A) , \tens\sigma\!{}_S \rangle\, dy - \dfrac 12 \int_Y \langle \tens C  \tens e(u_A), \tens e(u_A)\rangle\, dy
$$
and since from Theorem \ref{thrm:form-stress} $\tens\sigma\!{}_S  = \tens C \tens e(u_A)$, the last expression is equal to
$$
\int_Y \langle \tens e(u_A) , \tens\sigma\!{}_S \rangle\, dy
- \dfrac 12 \int_Y \langle \tens e(u_A) , \tens\sigma\!{}_S \rangle\, dy =
 \dfrac 12 \int_Y \langle \tens e(u_A) , \tens\sigma\!{}_S \rangle\, dy\, .
$$
Moreover, from Lemma \ref{thrm:dual-problem} below one gets
$$
 G(\tens\sigma\!{}_S) = g(\tens\sigma\!{}_S) = \dfrac 12 \int_Y \langle \tens D \tens\sigma\!{}_S , \tens\sigma\!{}_S \rangle = 
 \dfrac 12 \int_Y \langle \tens e(u_A) , \tens\sigma\!{}_S \rangle .
$$
Hence the conclusion of the Theorem follows.
\end{proof}

\begin{lemma}
\label{thrm:dual-problem}
For the function $G^\star : \LL^2_{\# s} \mapsto \RR\cup \{-\infty \}$ defined in (\ref{eq:G-star}) by 
$$
G^\star (\tens e) = \inf_{\tens\sigma \in \LL^2_{\# s}} \bigl\{ \int_Y \langle \tens e , \tens\sigma \rangle\, dy
+ \I_S(\Lambda \tens\sigma) \bigr\} - g^\star(\tens e) ,
$$ the dual problem (\ref{eq:dual_problem}) above
$$
 \displaystyle{\sup_{\tens e \in \LL^2_{\# s}} G^\star (\tens e)}
$$
may be written as 
$$
\displaystyle{\sup_{\tens e \in \LL^2_{\# s}} G^\star (\tens e)}= -  \displaystyle{\inf_{\tens e \in \TT(0)^\perp} k(\tens e)} ,
$$
where the function $k : \LL^2_{\# s} \mapsto \RR$ is defined in Theorem \ref{thrm:form-strain} by  
$
\displaystyle{k(\tens e) = \dfrac{1}{2} \dint_Y \langle \tens C \tens e, \tens e \rangle\, dy
  - \langle S, \dint_Y \tens e\, dy \rangle .}
$

\noindent Moreover, in the primal problem (\ref{eq:primal-problem-G}), the minimum value is equal to the maximum value in the dual problem and
$$
G(\tens\sigma\!{}_S)= \displaystyle{\inf_{\tens\sigma \in \LL^2_{\# s}} G(\tens\sigma)} = \displaystyle{\sup_{\tens e \in \LL^2_{\# s}} G^\star (\tens e)}
= G^\star(\tens e (u_A) ).
$$
\end{lemma}

\begin{proof}
The function $G^\star$ may be written as follows 
$$
G^\star (\tens e) = \inf_{\tens\sigma \in \LL^2_{\# s}} \bigl\{ \int_Y \langle \tens e , \tens\sigma \rangle\, dy
+ I_{\TT(S) } (\tens\sigma)\bigr\} - g^\star(\tens e) = $$
$$
\inf_{\tens\sigma \in \TT(S)} \bigl\{ \int_Y \langle \tens e , \tens\sigma \rangle\, dy \bigr\} - g^\star(\tens e) .$$
For $ \tens\sigma \in \TT(S) $
$$
 G^\star (\tens e)  = \displaystyle{\left\{\begin{matrix}
-\infty , \quad \mbox{ if } \tens e \notin \TT(0)^\perp , \\
\dint_Y \langle \tens e , S \rangle\, dy
- \dfrac 12 \int_Y \langle \tens C  \tens e , \tens e \rangle\, dy ,
\quad \mbox{ if } \tens e \in \TT(0)^\perp  ,
\end{matrix} \right.
}
$$
and therefore 
$$
\displaystyle{\sup_{\tens e \in \LL^2_{\# s}} G^\star (\tens e)}= \displaystyle{\sup_{\tens e \in \TT(0)^\perp } G^\star (\tens e)}
= \displaystyle{\sup_{\tens e \in \TT(0)^\perp } (-k (\tens e))} =  - \displaystyle{\inf_{\tens e \in \TT(0)^\perp } k (\tens e)} .
$$
On the other hand
$$
\inf_{\tens\sigma \in \LL^2_{\# s}} G (\tens\sigma) = 
\inf_{\tens\sigma \in \LL^2_{\# s}} \{ \dfrac 12 \dint \langle \tens D \tens\sigma , \tens\sigma \rangle\, dy, \}
$$
according to Theorem \ref{thrm:form-stress}
$$
=  \dfrac 12 \dint \langle \tens D \tens\sigma\!{}_S , \tens\sigma\!{}_S \rangle = G(\tens\sigma\!{}_S) 
$$
and, in view of Theorem \ref{thrm:form-strain}  and of the variational formulation (\ref{eq:var-form-displ-macrostress-mean}),
$$
= \dfrac 12 \int_Y \langle \tens C  \tens e(u_A), \tens e(u_A)\rangle\, dy =
\dfrac 12  \langle S, \int_Y \tens e(u_A)\, dy \rangle
- \dfrac 12 \int_Y \langle \tens C  \tens e(u_A), \tens e(u_A)\rangle\, dy
= -k(\tens e (u_A) ) = 
$$
$$
- \inf_{\tens e \in \TT^\perp_0} k (\tens e) = \sup_{\tens e \in \TT^\perp_0} G^\star (\tens e)=
G^\star(\tens e (u_A) ) = \sup_{\tens e \in \LL^2_{\# s}} G^\star (\tens e),
$$
which concludes the proof. 
\end{proof}

The minimization problem in stress (\ref{eq:primal-problem-G}) (primal problem) has also a dual problem in displacement. 
A Lagrangian associated to the dual problem in displacement is introduced in the next theorem.
The Lagrangian has a saddle-point that corresponds to the solutions of the primal and dual problems,
(\ref{eq:primal-problem-G}) and (\ref{eq:dual_problem}) respectively.

\begin{theorem}
\label{thrm:Lagrangean-stress-displ}
Consider the Lagrangian $\tilde \L : \LL^2_{\# s} \times LP_0 \mapsto \RR \cup \{ + \infty \}$ defined by
\begin{equation}
\tilde \L (\tens\sigma, v) = \dfrac 1 2 \int_Y \langle \tens D \tens\sigma , \tens\sigma  \rangle\, dy +
 \langle \Lambda \tens\sigma , v  \rangle
 - \I_S^\star (v) \quad \forall (\tens\sigma , v) 
\in \LL^2_{\# s} \times LP_0 .
\end{equation}
Then
$$
\displaystyle{\inf_{\tens\sigma \in \LL^2_{\# s}} G(\tens\sigma)=\inf_{\tens\sigma \in \LL^2_{\# s}}\sup_{v\in LP_0} \tilde\L (\tens\sigma, v) = 
\tilde \L (\tens\sigma\!{}_S, u_A) =
\sup_{v\in LP_0}\inf_{\tens\sigma \in \LL^2_{\# s}} \tilde \L (\tens\sigma, v)},
$$
where $\tens\sigma\!{}_S$ is the solution of both minimization problems in Theorem \ref{thrm:form-stress} and Theorem \ref{thrm:minimization.in.stress}  
and $\tens e (u_A)= \nabla_s u_A$ is the solution 
in Theorem \ref{thrm:form-strain}, provided $ S = \tens C^H A $.
\end{theorem}

\begin{proof}
For a given ${\tens\sigma \in \LL^2_{\# s}}$ one has
$$
\sup_{v\in LP_0} \tilde\L (\tens\sigma, v) = \dfrac 1 2 \int_Y \langle \tens D \tens\sigma , \tens\sigma  \rangle\, dy + 
\sup_{v\in LP_0} \{ \langle \Lambda \tens\sigma , v  \rangle
 - \I_S^\star (v) \}
$$
and using Theorem \ref{thrm:L-F-transf}
$$ = \dfrac 1 2 \int_Y \langle \tens D \tens\sigma , \tens\sigma  \rangle\, dy + 
 \sup_{v\in LP_0} \{  \int_Y \langle \tens\sigma - S , \nabla_s v \rangle\, dy =
 \dfrac 1 2 \int_Y \langle \tens D \tens\sigma , \tens\sigma  \rangle\, dy + I_{\TT(S)}(\tens\sigma ) = G(\tens\sigma ).
$$
Therefore 
$$
G(\tens\sigma )= \sup_{v\in LP_0} \tilde\L (\tens\sigma, v) \quad \hbox{ and consequently }
\displaystyle{\inf_{\tens\sigma \in \LL^2_{\# s}} G(\tens\sigma)=\inf_{\tens\sigma \in \LL^2_{\# s}}\sup_{v\in LP_0} \tilde\L (\tens\sigma, v)} .
$$
The infimum is attained in $\tens\sigma\!{}_S$\,:
$$
G(\tens\sigma\!{}_S) = \displaystyle{\inf_{\tens\sigma \in \LL^2_{\# s}} G(\tens\sigma) . }
$$
On the other hand 
$$
\tilde \L (\tens\sigma\!{}_S , u_A) = \dfrac 1 2 \int_Y \langle \tens D \tens\sigma\!{}_S , \tens\sigma\!{}_S  \rangle\, dy +
 \langle \Lambda \tens\sigma\!{}_S , u_A  \rangle
 - \I_S^\star (u_A) = \dfrac 1 2 \int_Y \langle \tens D \tens\sigma\!{}_S , \tens\sigma\!{}_S  \rangle\, dy +
 \int_Y \langle \tens\sigma\!{}_S - S , \nabla_s u_A \rangle\, dy 
$$
and using Theorem \ref{thrm:Donati.3} the last term vanishes, so
$$
= \dfrac 1 2 \int_Y \langle \tens D \tens\sigma\!{}_S , \tens\sigma\!{}_S  \rangle\, dy = G(\tens\sigma\!{}_S)
$$
and the inf-sup is attained in $(\tens\sigma\!{}_S , u_A) $.

Defining $\tilde G^\star : LP_0 \mapsto \RR\cup \{-\infty \}$ by 
$$
\tilde G^\star (v) = \inf_{\tens\sigma \in \LL^2_{\# s}} \bigl\{ \dfrac 1 2 \int_Y \langle \tens D\tens\sigma ,\tens\sigma \rangle\, dy +
\langle \Lambda \tens\sigma , v  \rangle \bigr\} - \dint_Y \langle S , \nabla_s v \rangle dy,
$$ 
the dual problem of (\ref{eq:primal-problem-G}) is
$$
 \displaystyle{\sup_{v\in LP_0} \tilde G^\star (v)} .
$$
Using the lemma below it turns out that
$$
\displaystyle{\sup_{v\in LP_0}\inf_{\tens\sigma \in \LL^2_{\# s}} \tilde \L (\tens\sigma, v)} = \sup_{v\in LP_0} \tilde G^\star (v) = 
\displaystyle{\inf_{\tens\sigma \in \LL^2_{\# s}} G(\tens\sigma)} = G(\tens\sigma\!{}_S),
$$
which concludes the proof.
\end{proof}

\begin{lemma}
\label{thrm:dual-problem-displ}
Let the function $\tilde G^\star : LP_0 \mapsto \RR\cup \{-\infty \}$ be defined by 
$$
\tilde G^\star (v) = \inf_{\tens\sigma \in \LL^2_{\# s}} \bigl\{ \dfrac 1 2 \int_Y \langle \tens D\tens\sigma ,\tens\sigma \rangle\, dy +
\langle \Lambda \tens\sigma , v  \rangle \bigr\} - \I_S^\star (v).
$$ 
Then the dual problem 
$$
 \displaystyle{\sup_{v\in LP_0} \tilde G^\star (v)}
$$
may be written as 
$$
\displaystyle{\sup_{v\in LP_0} \tilde G^\star (v)}= -  \displaystyle{\inf_{v \in LP_0} K(v)} ,
$$
where the function $K : LP_0 \mapsto \RR$ is defined in Theorem \ref{thrm:displ.2} by  
$
K(u) = \dfrac{1}{2} \dint_Y \langle \tens C \nabla_s u ,\nabla_s u \rangle\, dy
- \dint_Y \langle S ,  \nabla_s u \rangle \, dy  .
$

\noindent Moreover, in the primal problem (\ref{eq:primal-problem-G}), the minimum value is equal to the maximum value in the dual problem and
$$
G(\tens\sigma\!{}_S)= \displaystyle{\inf_{\tens\sigma \in \LL^2_{\# s}} G(\tens\sigma)} = \displaystyle{\sup_{v \in LP_0} \tilde G^\star (v)}
=  \tilde G^\star(- u_A ).
$$
\end{lemma}

\begin{proof}
For any $v\in LP_0$
$$
\tilde G^\star (v) = \inf_{\tens\sigma \in \LL^2_{\# s}} \bigl\{ \dfrac 1 2 \int_Y \langle \tens D\tens\sigma ,\tens\sigma \rangle\, dy +
\langle \Lambda \tens\sigma , v  \rangle \bigr\} - \I_S^\star (v) 
$$ 
and since from Theorem \ref{thrm:L-F-transf}, $\displaystyle{\I_S^\star (v) = \eta_S(v) = \dint_Y \langle S , \nabla_s v \rangle dy} $,
$$
\tilde G^\star (v)  = \inf_{\tens\sigma \in \LL^2_{\# s}} \bigl\{ \dfrac 1 2 \int_Y \langle \tens D\tens\sigma ,\tens\sigma \rangle\, dy +
\int_Y  \langle \tens\sigma , \nabla_s v \rangle \bigr\} - \dint_Y \langle S , \nabla_s v \rangle dy .
$$
The infimum is attained in $\tens\sigma = - \tens C \nabla_s v$ which implies that $ \tens D\tens\sigma = - \nabla_s v$, hence
$$
\tilde G^\star (v) = \bigl\{ \dfrac 1 2 \int_Y \langle \nabla_s v, \tens C \nabla_s v \rangle\, dy -
\int_Y  \langle \tens C \nabla_s v , \nabla_s v \rangle\, dy  \bigr\} - \dint_Y \langle S , \nabla_s v \rangle dy $$
$$=
-\dfrac 1 2 \int_Y \langle \nabla_s v, \tens C \nabla_s v \rangle\, dy -\dint_Y \langle S , \nabla_s v \rangle dy .
$$
Consequently
$$
\displaystyle{\sup_{v \in LP_0} \tilde G^\star (v)}= \displaystyle{\sup_{v \in LP_0} \tilde G^\star (-v)} =
\displaystyle{\sup_{v \in LP_0}}  \bigl\{-\dfrac 1 2 \int_Y \langle \nabla_s v, \tens C \nabla_s v \rangle\, dy 
+\dint_Y \langle S , \nabla_s v \rangle dy \bigr\} = - \displaystyle{\inf_{v \in LP_0} K(v) }
$$
On the other hand
$$
\inf_{\tens\sigma \in \LL^2_{\# s}} G (\tens\sigma) = 
\inf_{\tens\sigma \in \LL^2_{\# s}} \{ \dfrac 12 \dint \langle \tens D \tens\sigma , \tens\sigma \rangle\, dy \}
$$
according to Theorem \ref{thrm:form-stress}
$$
=  \dfrac 12 \dint \langle \tens D \tens\sigma\!{}_S , \tens\sigma\!{}_S \rangle = G(\tens\sigma\!{}_S) 
$$
and in view of Theorem \ref{thrm:displ.2} 2)  and of the variational formulation (\ref{eq:var-form-displ-macrostress})
$$
= \dfrac 12 \int_Y \langle \tens C  \nabla_s u_A, \nabla_s u_A \rangle\, dy =
 \langle S, \int_Y \nabla_s u_A \, dy \rangle
- \dfrac 12 \int_Y \langle \tens C  \nabla_s u_A , \nabla_s u_A\rangle\, dy
= -K( u_A ) $$
$$=  \tilde G^\star (- u_A) = - \displaystyle{\inf_{v\in LP_0} K(v) } = \displaystyle{\sup_{v\in LP_0} \tilde G^\star(v)}.
$$
\end{proof}

\section{Conclusions}
\label{sec:conclusions}

Several variational formulations of the cellular problem are presented, having as unknown
the displacement, the strain or the stress.
Although these formulations are equivalent, the formulation in stress is special in the sense that
the minimization problem in displacement and the minimization problem in strain may be viewed as
dual problems for the minimization problem in stress.
Thus, it is natural that in Theorems \ref{thrm:displ.2} 2) (formulation in displacement) and
\ref{thrm:form-strain} 2) (formulation in strain)
the functionals to minimize involve the macroscopic stress $S$ as datum.

Two lagrangians are constructed based on these primal-dual problems,
providing a stress-displacement approach and a stress-strain approach.
Each one of these approaches gives a more complete information about the cellular problem than
the stand-alone variational formulations, and are to be used according to the sought application.

Uzawa's method may be used to obtain the solution of the saddle point problem.
It uses a coupled iterative scheme along alternated directions.
Other, more performant, schemes, e.g.\ Arrow-Hurwicz, can also be used,
see \cite{AHU1958}. For numerical solution of saddle point problems see also \cite{BenGoLie2005}.

\section*{Aknowledgements}

This work is supported by national funding from FCT - Foundation for Science and Technology (Portugal),
project UIDB/04561/2020.

\bibliography{biblio}

\end{document}